\documentclass[10pt,a4paper]{amsart}
\usepackage{pictex}
\input amssym.def
\input amssym.tex
\title[Normal subgroups and congruence subgroups ]
{Identifying normal and congruence subgroups
}

\author{\SMALL Cheng Lien Lang}
\author{ \SMALL  mong lung lang}

\begin{document}

\baselineskip=12pt


\maketitle

\vspace{-0.3in}

\begin{abstract}
We study the arithmetic properties of the maps associated with
 the subgroups of finite indices  of the Hecke  group $G_q$.
  A necessary and sufficient condition for $X\subseteq G_q$
   to be normal is obtained. Such condition can be checked easily
    by {\bf GAP}.
    As a byproduct, one may determine the normaliser of $X$ in $PSL(2, \Bbb R)$
 when $q\ne 3,4,6$ and whether $X$ is congruence when $q\le 6$.


\end{abstract}

\section{Introduction}
\subsection{}Normal subgroups of the Hecke groups. Let $q\ge 3$ be a a fixed integer.
The (inhomogeneous) Hecke group
 $G_q$ is defined to be the maximal discrete subgroup
of  $PSL (2,\Bbb R)$
 generated by  $S$ and $T$, where $\lambda _q =2$cos$\,(\pi/q)$,
 $$ S =
\left (
\begin{array}{rr}
0 & 1 \\
-1 & 0 \\
\end{array}
\right ) \,,\,\,
  T = \left (
\begin{array}{rc}
1 & \lambda  _ q\\
0 & 1 \\
\end{array}
\right ) \,.\eqno(1.1)
$$
Normal subgroups  $X$ of finite indices of $G_q$ have been studied extensively,
 and various necessary and sufficient conditions have been imposed on $X$ for being normal
 (see [CS], [JS],  [Ma], [Mc], [N]
   for examples).
   To the best of our knowledge,  the verification of those conditions is quite
   delicate and can be lengthy.
       For instance, while it is well known that $X$ is normal if and only if
        the action of $G_q$ on $G_q/X$ is fixed point free,
        the actual checking is very tedious as it involves
         some lengthy membership test.


\subsection{}The main results.
 The main purpose of this article is to associate with each subgroup $X \subseteq G_q$
  (given by generators and special polygons as in subsection 1.3)
   a  finite  subgroup $G_{X}$ in the sense of
    Millington [Mi1, Mi2] and Atkin and Swinnerton-Dyer [AS] such that
  \begin{enumerate}
  \item[(i)]  the generators of $G_X = \left <r_1, r_2\right >$ can be determined
   directly from the special polygon of $X$ without matrix multiplication
     or any  membership test and that
   the order of $G_{X}$ can be calculated by {\bf GAP},
   \item[(ii)]  whether
    $X$ is normal in $G_q$ can be  determined by whether the order of $G_{X}$ meets
     the index $[G_q : X]$ (see Proposition 6.2).
     \end{enumerate}
 As a byproduct, one may determine the normaliser of $X$ in $PSL(2, \Bbb R)$
when $q\ne 3, 4,6$ (see subsection 7.2)
and whether $X$ is congruence when $q\le 6$ (see Section 8).
 Note that the normalisers of
  subgroups of finite indices of $PSL(2, \Bbb Z)$
  in $PSL(2, \Bbb R)$ can be determined by  studying the {\em big picture} of Conway
   (see [L]). The normalisers of subgroups of finite indices of $G_4$ and
   $G_6$ in $PSL(2, \Bbb R)$ are open to us. 
   The permutations $\sigma_0$ and
    $\sigma_1$ of He and Read 
    in their study of {\em dessins d'enfants} (pp.12 of [HR]) can be calculated easily by our results
     as $\sigma_0 = r_1$ and $\sigma_1= r_2$.

\subsection{} The main assumption.
The main assumption of this article is that subgroups of finite indices of $G_q$
 are described
 as in  Kulkarni [K]. Equivalently,
 in terms of fundamental domains and sets of independent generators in the sense of
  Rademacher.
The existence of such description is
  guaranteed
  by Kurosh Theorem [Ku]. To be more precise, it is  proved in [LLT1] that every subgroup
$X$   of finite
  index of $G_q$ has a fundamental domain $M_X$ which is a special polygon and that the set
   of side
   pairings of
   the Hecke-Farey symbols associated with
   the special polygon $M_X$ is a set of independent generators of $X$.

  \subsection{} A simple  idea. Jones and Singerman in their beautiful paper [JS]
   give a close relationship
   between  subgroups  of the {\em triangle  groups}  and
    {\em maps}.
    Following their insight,
     we make the special polygon $M_X$ into a  map and investigate the
      monodromy
       group $G_X=\left <r_1, r_2\right >$ of $M_X$. As a homomorphic
        image of $G_q$,  $G_X$
        acts on the set of cosets $G_q/G_0$, where
        $G_0$ is a map subgroup (see Definition 5.1).
           The key of our study is that  $M_X$ is constructed
        in such a way that
          (i) $G_0$ is a conjugate of $X$,
         (ii)
         the permutations  $r_1$ and $r_2$  are compatible with the
   action of $S$ and $ R = ST^{-1} $ on $G_q/X$ (see Lemmas 5.2-5.4),
    (iii) $r_1$ and $r_2$ can be written down without any membership test
     (see Examples 7.1-7.5).
    As a consequence,
      the conditions (i) and (ii) of subsection 1.2 are
     satisfied and we are able to determine the normality of $X$
      by checking whether
       $|G_X|/[G_q : X]$ is one (Proposition 6.2).

\subsection{} The organisation.
 In Section 2, we give  a lemma that will be used in Section 6. See [W] for more detail.
  The main theme of Section 3 is the geometrical construction of
   a special polygon (map) and a set of independent generators of $X$.
      Section 4 gives a very short study of the basic properties
       of maps that is useful for our study. Section 5
       studies the
        action of $G_X$ on $G_q/X$.
        Section 6 gives a simple test that enables us to identify normal subgroups
         of $G_q$.
        Section 7 demonstrates how normal subgroups can be identified and
         how the normaliser
        $N(X) $  can be calculated.
      Example 7.5 gives a normal non-congruence subgroup of $PSL(2, \Bbb Z)$
    of index 42 which is the smallest possible index of a normal non-congruence
     subgroup of $PSL(2, \Bbb Z)$.
           Section 8  studies the congruence
          subgroup problem and
          Appendix A studies the torsion normal subroups.

\smallskip
It is worthwhile to note that
 $ G_0(2)=\{
       (x_{ij} )\in G_6 : x_{12} \equiv 0\,(mod\,\,2)\}$
        has index 3 in $G_6$ (see (ii) of Example 7.4). This is a little
         surprise to us as one  expects,
           in
         line with the modular group case where a special polygon of $\Gamma_0(2)$
          is a 3-gon, that
     a special polygon of $G_0(2)$
          is   a 6-gon
              and that
             $[G_6 : G_0(2)] = 6$.




\section {a Technical Lemma  }
Throughout the section, $\Omega$ is a finite set and $S_{\Omega}$ is the
 symmetric group on $\Omega$.




\smallskip
\noindent {\bf Lemma 2.1.} {\em
Suppose that $G$ acts  transitively on $\Omega$.
      Then $C_{S_{\Omega}}(G) \cong N_G(G_d)/G_d$,
       where $G_d$ is the one point stabiliser of $d \in \Omega$.}

       \smallskip
       \noindent {\em Proof.} Since $G$ is transitive,
       The action of $G $ on $\Omega$ is isomorphic to the   action of $G$ on the set of
        cosets $G/G_d$. Without loss of generality, we may assume that
        $\Omega = G/G_d$ and that $x(G_d) = xG_d$ for $x \in G$.
            Let $x\in C_{S_{\Omega}}(G)$. Then $x(G_d) = e_xG_d\in G/G_d$ for some $e_x \in G$.
     For each $g \in G_d$, one has $gx(G_d) = xg(G_d) = x(G_d)$.
     Hence $ge_xG_d = e_x G_d$. This implies that
     $e_xG_d \in N_G(G_d)/G_d$. As a consequence, one can show that
     $C_{S_{\Omega}}(G) \cong N_G(G_d)/G_d$ by studying the
      homomorphism  $\Phi\, :\, C_{S_{\Omega}}(G) \to N_G(G_d)/G_d$ defined
       by $\Phi(x) = e_x^{-1} G_d$. Note that for each $  r \in N_G(G_d)$,
        the permutation defined by $\sigma (g G_d) = gr^{-1} G_d$ commutes with $G$
         which implies that $\Phi$ is surjective  ($e_{\sigma} = r^{-1}$, $\Phi (\sigma) = r^{}G_d$).
    \qed


\smallskip
\noindent {\bf Corollary 2.2.} {\em Suppose that $G$ acts transitively and freely
$($equivalently, $G_{x} = 1$ for
    all $x \in \Omega)$
on
 $\Omega$. Then $G \cong C_{S_{\Omega}}(G)$.}


\section{Fundamental domains of subgroups of $G_q$  }

\subsection{}In [K], Kulkarni applied a combination of geometric
 and arithmetic methods to show that one can produce
 a set of independent generators in the sense of
 Rademacher for the congruence subgroups of the modular group,
 in fact for all subgroups of finite indices.
 His method can be generalised to all subgroups of
 finite indices of the Hecke groups $G_q$.
  In short, for each subgroup $X$ of finite index of $G_q$,
one can associate with  $X$
a set of Hecke-Farey symbols (HFS)
$\{-\infty, x_0, x_1, \cdots, x_n, \infty\},$
 a special  polygon (fundamental domain) $M_X$,
 and  an additional structure
 on each consecutive pair of $x_i$'s of the four types described
 below :
$$ {x_i} _{_{_{\smile}}} \  \hspace{-.37cm} _{_{_{_{_{_{_{\circ}}}}}}}
  x_{i+1} ,\,\,
{x_i} _{_{_{\smile}}} \  \hspace{-.37cm} _{_{_{_{_{_{_{\bullet}}}}}}}
 x_{i+1} ,\,\,
  {x_i} _{_{_{\smile}} }\  \hspace{-.37cm} _{_{_{_{_{_{_{a}}}}}}}  x_{i+1},\,\,
  {x_i}\,\, _{_{_{\smile}} }   \hspace{-.24cm} _{_{_{_{_{_{_{e_{r}}}}}}}}  x_{i+1}.
  \eqno(3.1)
 $$
where $a$ is a nature  number and $1 <r <q$ is a divisor of $q$. Each nature number $a$ occurs
 exactly twice or not at all.
Similar to the modular group,
 the actual values of the $a$'s is unimportant: it is
 the pairing induced on the consecutive pairs that matters.

 \begin{enumerate}

\item[(0)]  A set of Hecke-Farey symbols (HFS) is a  finite sequence of cyclically arranged
  numbers
 $\{-\infty, x_0, x_1, \cdots, x_n, \infty\}$
  such that
  \begin{enumerate}
  \item[(a)] $ x_i \in \Bbb Q (\lambda_q)$, $x_i=0$ for some $i$, $0\le i \le n$,
  \item[(b)]   $x_i = a_i/b_i$ is in reduced form  for every $i$ (letting $x_{-1} =-\infty
   = -1/0$ and $x _{n+1} = \infty = 1/0$).
  \end{enumerate}

  \item[(i)] The side pairing  $\circ$ is the following  elliptic element of order 2 that pairs
  the even line $(a/b, c/d)$ with itself ($cb-ad=1$). The trace of such an element is 0.
   Recall that the $G_q$-translates of the hyperbolic line $(0, \infty)
  = (0, i)\cup (i,\infty)$ are called
   the even lines.
  {\small   $$
\left (
\begin{array}{rr}
c & a \\
d & b \\
\end{array}
\right )
  \left (
\begin{array}{rr}
0 & 1\\
-1 & 0 \\
\end{array}
\right )
\left (
\begin{array}{rr}
c & a \\
d & b \\
\end{array}
\right ) ^{-1}
\in G_q \,.\eqno(3.2)
$$ }
   \item[(ii)]
  The side pairing  $\bullet$ is the following  elliptic element of order $q$ that pairs
  the odd line $(a/b, c/d)$ with itself ($cb-ad=1$).
   The absolute value of the trace of such an element is
  $\lambda_q$.
   Recall that  the $G_q$-translates of the hyperbolic line $(0, e^{\pi i/q})
   \cup ( e^{\pi i/q} , \infty)$ are called
   the odd lines.
     {\small   $$
\left (
\begin{array}{rr}
c & a \\
d & b \\
\end{array}
\right )
  \left (
\begin{array}{rc}
0 & 1\\
-1 & \lambda_q \\
\end{array}
\right )
\left (
\begin{array}{rr}
c & a \\
d & b \\
\end{array}
\right ) ^{-1}
\in G_q \,.\eqno(3.3)
$$ }
     \item[(iii)] The two even lines $(a/b, c/d)$ and $(u/v, x/y)$ with the label $a$
   are paired together by the following  element of infinite order ($cb-ad=1$,
   $vx-yu=1$).
{\small   $$ \hspace{1.5cm}
\left (
\begin{array}{rr}
u & -x \\
v & -y \\
\end{array}
\right )
  \left (
\begin{array}{rr}
c & a\\
d & b \\
\end{array}
\right )^{-1}\in G_q \,.\eqno(3.4)
$$ }
\item[(iv)] The side pairing $e_r$ pairs $x_i= u/v$ and $x_{i+1}=x/y$
($vx-yu\ne1$)  as follows.
 There exists $y_2 < y_3 < \cdots < y_{q-q/r}$ such that
  $x_i, y_2, \cdots, y_{q-q/r}, x_{i+1}$ are consecutive entries of a $q$-gon $P$.
  Let $a/b$ and $c/d$  (in reduced forms)  be the smallest and largest entries of $P$.
  Then $x_i$ and $x_{i+1}$ are paired by
 {\small   $$
\left (
\begin{array}{rr}
c & a \\
d & b \\
\end{array}
\right )
  \left (
\begin{array}{rc}
0 & 1\\
-1 & \lambda_q \\
\end{array}
\right )^{q/r}
\left (
\begin{array}{rr}
c & a \\
d & b \\
\end{array}
\right ) ^{-1}
\in G_q \,.\eqno(3.5)
$$ }
 See subsection 3.3 and Example 7.4 for such side parings of order $r$.
 \item[(v)]
The special polygon $M_X$ associated with  the HFS is a fundamental
 domain of $X$.   It is a union of $q$-gons and $s$-gons
  (each $\bullet$ gives a 1-gon and each $e_r$ gives a $q/r$-gon, see Example 7.4).
  Each $s$-gon has $s$ special triangles.   The side pairings
 $ I= \{ \sigma_1,  \sigma_2 \cdots ,  \sigma_m \}$
associated with the  HFS
 is a set of independent generators of $X$.
The number $d$ of  special triangles
 (a special triangle is  a fundamental domain of $G_q$) of the special polygon is  the index
 of the
   subgroup $X$.
   \item[(vi)] The set of independent generators consists of
  $s$ matrices of infinite order, where $s$ is the number of the nature number $a$'s in the
   Hecke-Farey smybols.
   \item[(vii)] The subgroup $X$ has
   $\tau _2$
   (the number of the circles $\circ$ in HFS) inequivalent classes of elliptic
    elements of order 2  that are conjugates of $S$. Each class has exactly one representative in $I$.
     \item[(viii)] The subgroup $X$ has
   $v_q$
   (the number of the bullets $\bullet$ in HFS) inequivalent classes of elliptic
    elements of order $q$ that are conjugates of $ST^{-1}$. Each class has exactly one representative in $I$.
     \item[(ix)] The subgroup $X$ has
   $v_r$
   (the number $e_r$'s in HFS) inequivalent classes of elliptic
    elements of order $r$ that are conjugates of $(ST^{-1})^{q/r}$. Each class has exactly one representative in $I$.
     \item[(x)]
      The Hecke-Farey symbols (cusps)
     can be partitioned into $v_{\infty}$ classes under the action of the
      set of independent generators,
       such equivalence classes are called  the {\em vertices} of $M_X$.
       It follows that  $M_X$ has $v_{\infty}$ vertices.
       The  equivalence classes of even lines are called the {\em edges} of $M_X$.
        The $s$-gons (see (v)) are called the {\em faces} of $M_X$.
        The set of vertices, edges and faces are denoted by $V$, $E$
         and $F$ respectively.

     \item[(xi)] The degree of the  vertex $x$,
      denoted by $w(x)$,  is the number of even lines
      in $M_X$ that comes into $x$.
       Algebraically, it is the smallest positive integer
       $m$ such that $\pm T_q^m$ is conjugate in $G_q$ to an element of $X$ fixing $x$.
        Denote by $V $ the set of  vertices  of $X$. Then
         $ [G_q : X] = \sum _{x\in V} w(x) $.
 \item[(xii)] The genus $g$ of $M_X$ is given by the following
   Riemann Hurwitz formula.
$$ 2g-2 + \tau_2/2 + \sum v_r(1-1/r) + v_{\infty} = [G_q : X] (1/2-1/q)
,\eqno(3.6)$$
where the sum is taken over the set of all positive divisors of $q$.

        \item[(xiii)] Let $P$ be an ideal $q$-gon with cusps $\{x_1, x_2, \cdots , x_q\}$
 (if $\infty$ is a cusp, then $x_1 = -\infty = -1/0$ if $P$ lies to the left
  of the $y$-axis and $x_q=\infty=1/0$ if $P$ lies to the right of the $y$-axis).
   Let $x_i= a_i/b_i$ be in reduced form. Then
   $$    \left (\begin{array}{cc}
                            a_1 & -a_q      \\
                            b_1 & -b_q\\
                             \end{array}\right ) A^i =
                         \left (\begin{array}{cc}
                               a_{i+1} &   a_i     \\
                               b_{i+1} &  b_i\\
                             \end{array}\right ),\,\,
                             \mbox{ where }
                     A =
       \left (\begin{array}{cc}
                            \lambda_q  & 1     \\
                                    -1 & 0\\
                             \end{array}\right )  . \eqno(3.7)$$

   \end{enumerate}
 In the case $q$ is a prime or $X$ is torsion free, the side pairing $e_r$ does not exist
 and the entries of the Hecke-Farey symbols in (0) of the above satisfy
   $a_{i+1}b_i - a_i b_{i+1} = 1$ for all $i$.
See [LLT1] for more detail.
The numbers in (vi)-(xii)  are called the
 {\em geometric invariance} of $X$.

\subsection {Subgroups of $PSL(2, \Bbb Z)$  of index 7}
  It is well known that $PSL(2, \Bbb Z)$ has
  altogether 42 subgroups of index 7  (see  [R], [LLT2]).
  Since $PSL(2, \Bbb Z_7)$ has 14 subgroups of index 7,
  fourteen  of the 42 subgroups of index 7  are congruence of level 7. The remaining 28
  are non-congruence.
 A special polygon of such a subgroup  $A$ consists of two  3-gons and one
 special triangle.
 Since a special polygon of $A$ has a special triangle,
  $A$ has an element $g$ of order 3. Suppose that $A$ is normal. It follows that
   $A$ contains all the conjugates of $g$. In particular, $ST^{-1}, T^{-1}S \in A$.
    Hence a special polygon of $A$ is a union of two special triangles
      $ \{  {-\infty }_{_{_{\smile}} }\  \hspace{-.37cm}_{ _{_{_{_{_{_{_{\bullet}}}}}}} }
 {0}
_{_{_{\smile}} }\  \hspace{-.37cm}_{ _{_{_{_{_{_{_{\bullet}}}}}}}}
\infty
\}$.    This contradicts the fact that $A$ has index 7. Hence $A$ is not normal.
 As a consequence, $A$ is self-normalised and  has 7 conjugates.
  We shall list for each conjugacy class a representative as follows.
     Note that the subscripts refer to the degrees of the vertices.
$$M_{7} = \{  {-\infty }_{_{_{\smile}} }\  \hspace{-.37cm}_{ _{_{_{_{_{_{_{\bullet}}}}}}} }
 {-1}
_{_{_{\smile}} }\  \hspace{-.37cm}_{ _{_{_{_{_{_{_{\circ}}}}}}}}
0
 _{_{_{\smile}}}\  \hspace{-.36cm} _{_{_{_{_{_{_{_{\circ}}}}}}}}
 {1}
_{_{_{\smile}} }\  \hspace{-.37cm}_{ _{_{_{_{_{_{_{\circ}}}}}}}}
 {\infty}
 \} .\eqno(3.8)$$
$$M_{7} '= \{  {-\infty }_{_{_{\smile}} }\  \hspace{-.37cm}_{ _{_{_{_{_{_{_{\circ}}}}}}} }
 {-1}
_{_{_{\smile}} }\  \hspace{-.37cm}_{ _{_{_{_{_{_{_{\bullet}}}}}}}}
0
 _{_{_{\smile}}}\  \hspace{-.36cm} _{_{_{_{_{_{_{_{\circ}}}}}}}}
 {1}
_{_{_{\smile}} }\  \hspace{-.37cm}_{ _{_{_{_{_{_{_{\circ}}}}}}}}
 {\infty}
 \}. \eqno(3.9)$$
Subgroups (3.8) and (3.9) are congruence subgroups of level 7.
 The following
 subgroups
 are non-congruence.
 $$M_{1,6} = \{  {-\infty }_{_{_{\smile}} }\  \hspace{-.37cm}_{ _{_{_{_{_{_{_{\bullet}}}}}}} }
 {-1}\,
_{_{_{\smile}} }\  \hspace{-.37cm}_{ _{_{_{_{_{_{_{\circ }}}}}}}}
0
 _{_{_{\smile}}  }\  \hspace{-.37cm} _{_{_{_{_{_{_{_{1}}}}}}}}
\,{1}
 _{_{_{\smile}}}\  \hspace{-.36cm} _{_{_{_{_{_{_{_{1}}}}}}}}
 {\infty}
  \}.
    \eqno(3.10)$$
 $$M_{1,6}' = \{  {-\infty }_{_{_{\smile}} }\  \hspace{-.37cm}_{ _{_{_{_{_{_{_{\circ}}}}}}} }
 {-1}\,
_{_{_{\smile}} }\  \hspace{-.37cm}_{ _{_{_{_{_{_{_{\bullet }}}}}}}}
0
 _{_{_{\smile}}  }\  \hspace{-.37cm} _{_{_{_{_{_{_{_{1}}}}}}}}
\,{1}
 _{_{_{\smile}}}\  \hspace{-.36cm} _{_{_{_{_{_{_{_{1}}}}}}}}
 {\infty}
  \}.
    \eqno(3.11)$$
       $$M_{3,4} = \{  {-\infty }_{_{_{\smile}} }\  \hspace{-.37cm}_{ _{_{_{_{_{_{_{\bullet}}}}}}} }
 {-1}
_{_{_{\smile}} }\  \hspace{-.37cm}_{ _{_{_{_{_{_{_{1}}}}}}}}
0
 _{_{_{\smile}}}\  \hspace{-.36cm} _{_{_{_{_{_{_{_{\circ}}}}}}}}
 {1}
_{_{_{\smile}} }\  \hspace{-.37cm}_{ _{_{_{_{_{_{_{1}}}}}}}}
 {\infty}\}.
 \eqno(3.12)$$
 $$
M_{2,5} = \{  {-\infty }_{_{_{\smile}} }\  \hspace{-.37cm}_{ _{_{_{_{_{_{_{\bullet}}}}}}} }
 {-1}
_{_{_{\smile}} }\  \hspace{-.37cm}_{ _{_{_{_{_{_{_{1}}}}}}}}
0
 _{_{_{\smile}}}\  \hspace{-.36cm} _{_{_{_{_{_{_{_{1}}}}}}}}
 {1}
_{_{_{\smile}} }\  \hspace{-.37cm}_{ _{_{_{_{_{_{_{\circ}}}}}}}}
 {\infty}
 \}. \eqno(3.13)$$
The geometric invariance  of the above subgroups can be determined by (vi)-(xii)
 of subsection 3.1. For instance, the invariance of the subgroup (3.8) are $\tau_2 = 3$, $v_3 =1$,
  $v_{\infty} = 1$, the degree of the only vertex is $7$. The genus of the surface
   is zero. $M_{1,6}$ will be used to construct a normal non-congruence subgroup of
    index 42 (see Example 7.5).

\subsection{ A subgroup of  the Hecke group $G_6$} Let $\lambda = \sqrt 3$ and let $M_X$ be given as follows.
$$
  M_X= \{
   {-\infty }_{_{_{\smile}} }\  \hspace{-.37cm}_{ _{_{_{_{_{_{_{\circ}}}}}}} }
 \,\,
   {0/1 }_{_{_{\smile}} }\  \hspace{-.37cm}_{ _{_{_{_{_{_{_{\circ}}}}}}} }
 {1/\lambda}\,\,
_{_{_{\smile}} }\  \hspace{-.32cm}_{ _{_{_{_{_{_{_{\circ}}}}}}}}
\lambda/2 \,\,
_{_{_{\smile}} }\  \hspace{-.32cm}_{ _{_{_{_{_{_{_{e_2}}}}}}}}
\infty
 \}. \eqno(3.14)$$

\noindent An easy calculation (see (xiii) of subsection 3.1) shows that
 $\{
 \lambda/2,2/\lambda,\lambda/1, \infty\}$ are consecutive  entries of the 6-gon
 $P=\{ -\infty, 0/1, 1/\lambda, \lambda/2,2/\lambda,\lambda/1, \infty\}$. By (i) and (iv) of subsection 3.1,
  a set of independent generators is given by
      $$
    X = \left < S,
\left (
\begin{array}{rr}
1 & 0 \\
 \lambda & 1 \\
\end{array}
\right )
 S
\left (
\begin{array}{rr}
1 & 0\\
\lambda & 1 \\
\end{array}
\right ) ^{-1},
\left (
\begin{array}{rr}
\lambda & 1   \\
2 & \lambda \\
\end{array}
\right ) S
\left (
\begin{array}{rr}
\lambda  & 1\\
2 & \lambda \\
\end{array}
\right ) ^{-1},
\left (
\begin{array}{rr}
0 & 1\\
-1 & \lambda \\
\end{array}
\right ) ^{3}\right >.\eqno(3.15)
$$
Note that $[G_6 : X] = 3$. See Example 7.4 (Figure 3) for a picture  of
$M_X$.

\section{Maps of subgroups of $G_q$}

\subsection{ $M_X$ is a map}  Let $M_X$ be a special polygon  of $X$ given as in Section 3.
$M_X$ is an orientation-preserving compact Reimann surface.  Recall  that
 (see (x) of subsection 3.1)
\begin{enumerate}
\item[(i)] $V = $ set of vertices of $M_X =$ the set of equivalence classes of cusps of $M_X$,
\item[(ii)] $E=$ set of edges of $M_X =$   the
  set of  equivalence classes of even lines of $M_X$,
  \item[(iii)] $F=$ set of faces of $M_X=$ the set of
   $r$-gons  of  $M_X$.

    \end{enumerate}

    An even line $(x,y)= (a/b, c/d)$ paired
   by an elliptic element of order 2 (see (i) of subsection 3.1) gives an edge with one vertex only.
    Such an edge is called a {\em free edge}.  A non-free edge $(x,y)$ with only one
     vertex is called a {\em loop}
      (the end points $x$ and $y$ are paired by a set of side
     pairings of orders $\ge 3$). An  edge  with two vertices
is  called a {\em segment}.

     It is clear
 that the faces are simply connected. As a consequence, $M_X =(V, E, F)$
  is a  {\em map} (a map on a compact orientable surface $\mathcal S$ is an embedding of a
   finite connected graph $\Gamma$ on $\mathcal S$ such that the connected components of $\mathcal S \setminus \Gamma$
    are simply connected,
  see [CS] for more detail). $\Gamma_X = (V, E)$ is called the graph of $M_X$.


\smallskip



\subsection{Directed graph and  darts of $M_X$}
Let $e$  be an edge of the map $M_X$.
We may associate with the  edge $e$ two  directed edges if $e$ is not a free edge (see
  subsection 4.1) and
 one directed edge if $e$ is a free edge as follow.
\begin{enumerate}
\item[(i)]
Suppose that $e$ is a loop.  Then one of the
 directed edge travels clockwise and the other counter-clockwise.
 In the case $e$ is a segment with vertices $x$ and $y$, one travels from $x$ to $y$ and
  the other travels from $y$ to $x$.
  \item[(ii)]
  Suppose that  $e$ is a free edge.
Then we associate with $e$ one directed edge only.  Such a directed  edge always points
  to the vertex of the edge.
  \item[(iii)]
 A directed edge is called a
{\em  dart}.
   Denote by $\Omega$ the set of all darts.

 \end{enumerate}


\subsection{The cosets $G_q/X$ and the darts $\Omega$ have the same cardinality}
  Let $\Phi$ be the special triangle
 with vertices  $\infty$, 0 and $e^{\pi i/q}$. Then $\Phi$ is a fundamental domain of $G_q$.
 Suppose that $G_q = \cup _{k=1}^r Xg_i$.
 Replace $g_k$ by $x_kg_k$ for some $x_k\in X$ if necessary,
  we may assume that
 $\cup_{k=1}^r  g_k\Phi  = M_X$ (as a special polygon for $X$).
  Each dart $\hat d_k \in \Omega$  belongs to a unique special triangle $g_k\Phi
   \in M_X$ (following the orientation).
   As a consequence, there exists a
     one to one correspondence $\sigma$ between $ \Omega $ and $G_q/X$
     defined by  $\sigma(\hat d_k) = Xg_k$.
Note that each non-free  edge belongs to exactly two special triangles.

\section{  Monodromy groups of subgroups of $G_q$ }

\subsection{Monodromy groups of $M_X$}
Let $\Omega$ be the set of darts of $M_X$. Define $r_1$ to be the permutation
 that transposes the two darts of each loop and segment, and fixes the single dart of
  each free edge.  Around each vertex $v$ of $M_X$, the orientation of $M_X$
  imposes a cyclic ordering of the darts pointing towards $v$ and $r_0$ is the permutation
   with these as their disjoint cycles.
   Define $r_2 = r_1^{-1} r_0^{-1}$. $r_2$ is a permutation whose cycles correspond
    to the faces of $M_X$, again following the orientation.
   Define
               $$G_{X}  = \left <r_1, r_2\right >  =   \left <r_0, r_1\right >
                .\eqno(5.1)$$

\noindent
 The permutations $r_1$ and $r_2$ can be written down easily from $M_X$.
See Examples 7.1 and 7.2 for the actual readings of $r_1$ and $r_2$.
  Note that  the   incidence relations of $(V, E, F)$ are completely determined by $r_0$ and $r_1$
 and that every dart appears in exactly one cycle of the $r_i$'s ($0\le i\le 2$).
Let $R = ST^{-1}$.  Since $\{S, R \}$ is a set of independent generators of $G_q$, $o(R) = q$,
 and the order of $r_2$ is a divisor of $q$, the map
 $f$ defined by $$f(S) = r_1,\,f(R) = r_2\eqno(5.2)$$
  is a surjective  homomorphism
  from $G_q$ to $G_X$.
  It is clear from (5.2) that $G_q$ acts on $\Omega$ and  that the action of $G_X$ on $\Omega$ is isomorphic to the action of
   $G_q$ on $G_q/G_0$, where $G_0$ is any one point stabiliser of the action of $G_q$
    on $\Omega$.

    \smallskip
    \noindent {\bf Definition 5.1.}
    The group $G_X$ defined in (5.1) is called
     the {\em monodromy group} of $X$.
     The one point stabilisers of the action of
     $G_q$ on $\Omega$ are called the {\em map subgroups} of  $M_X$.
   Denote by $Y$ the kernel of $f$. $Y$ is called the {\em normal
   map subgroup} associated with $M_X$.

\subsection{Map subgroups of $X$}

Since the map $M_X$ is a special polygon of $X$,
 it is very natural to expect that the map subgroups of $M_X$ are the
  conjugates of $X$ (Lemma 5.4) and that
   the action of $S$ and $R$ on $G_q/X$ can be described by the action
    of $G_X$ on $\Omega$ (Lemmas 5.2 and 5.3).

 \smallskip
 \noindent {\bf Lemma 5.2.} {\em Let $\hat d_k, g_k$ and $ r_1$ be given as in subsections $4.3$
  and $5.1$.
 Suppose that $  r_1 \hat d_i  =f(S) \hat d_i = \hat d_j$. Then $Xg_i S = X g_j$.}

\smallskip
\noindent {\em Proof.}
Following the definition of $f(S)$, $ \hat d_i $ and $  \hat d_j$ are the directed edges associated with
  the same edge  $e$ and that $\hat d_i \in g_i\Phi$, $\hat d_j \in g_j\Phi$ (see subsection 4.3 for notation).

 \smallskip
Since $\hat d_i \in g_i \Phi$, one has
 $e \in g_i\Phi$.
 Hence $g_i^{-1} e \in \Phi
\cap S\Phi$.
Equivalently,  $ e \in g_i \Phi
\cap g_iS\Phi$.
This implies that $\hat d_i \in g_i\Phi$ and that
 $\hat d_j \in g_iS\Phi$.
 Since each dart belongs to a unique special triangle of $M_X$ (see subsection 4.3),
 one has $  x_ig_i S\Phi = g_j\Phi $ for some $x_i \in X$.
  Since $G_q$ acts freely on the $G_q$-translates of $\Phi$, one has $x_ig_i S =g_j$.
   Hence $Xg_i S = Xg_j$.
\qed

\smallskip

\noindent {\bf Lemma 5.3.} {\em Let $\hat d_k, g_k$ and $ r_2 $ be given as in subsections $4.3$
 and $5.1$.
Suppose that $ r_2 \hat d_i =f(R) \hat d_i = \hat d_j$. Then $Xg_i R = X g_j$.}

\smallskip
\noindent {\em Proof.}
Note first that  $\hat d_i \in g_i \Phi$,
 $\hat d_j \in g_j \Phi$  and that
 $ \hat d_i $ and $  \hat d_j$ (following the orientation) are consecutive terms of an $s$-gon of $M_X$.
  Since  $\hat d_i \in g_i \Phi$,
one has   $g_i^{-1} \hat d_i\in\Phi$. Hence
 $g_i^{-1}\hat  d_i = \hat e$, where $e$ is the even line $(0,\infty) \in \Phi$
  and $\hat e$ is a directed edge associated with $e$.
  It follows that
    $g_i^{-1} \hat d_i = \hat e$ is a term of an $s$-gon $F_s$ of $g_i ^{-1}Xg_i$.
     The term next to  $g_i^{-1} \hat d_i = \hat e$ (following the orientation)
      in $F_s$ is $R \hat e$.
       This implies that $ g_iR \hat e = \hat d_j  \in g_j \Phi$.
        As a consequence, one has
        $$\hat d_j  =g_iR \hat e =
        g_iR g_i^{-1} g_i \hat e = g_iR g_i^{-1} \hat d_i \in g_iR g_i^{-1} g_i\Phi  =
         g_iR \Phi   .\eqno(5.3)$$
Since each dart belongs to a unique special triangle of $M_X$ (see subsection 4.3),
 one has $x_ig_i R\Phi =  g_j\Phi $ for some $x_i \in X$. Hence $Xg_i R = Xg_j$.\qed


 \smallskip
 \noindent
{\bf Lemma 5.4.} {\em The map subgroups of $M_X$ are conjugates of $X$.
The action of
    $G_X$ on $\Omega$ is isomorphic to the action of $G_q$
     on $G_q/X$. Further, the normal map subgroup $Y=\cap\, gXg^{-1}$ is  a subgroup of $X$. }

\smallskip
\noindent {\em Proof.} $G_q$ acts on $\Omega$ by (5.2). Let $G_0$ be the stabiliser of
 $\hat d_i \in \Omega$. It follows easily that the action of $G_q$ on $G_q/G_0$ is isomorphic to the action of $G_X$ on $\Omega$. In particular,
   one has $|G_q/G_0| =|\Omega|$.

\smallskip
 Suppose that $f(g) \hat d_i = \hat d_i$.
  By Lemmas 5.2 and 5.3,
one has  $Xg_i g = Xg_i$.
  Equivalently, $ g \in g_i^{-1} Xg_i$. This implies
   that $G_0$ (the stabiliser of $\hat d_i$) is a subgroup of
     $ g_i^{-1} Xg_i$. Since $ G_0 \subseteq  g_i ^{-1}Xg_i$,
     $|G_q/G_0| =|\Omega|$
      and
      $|G_q/g_i^{-1} Xg_i|= |G_q/X| = |\Omega|$ (see subsection 4.3), one concludes
       that $G_0 =g_i^{-1} Xg_i$.
    As a consequence,  the map subgroups (see Definition 5.1) are conjugates of $X$ and
     the action of $G_X$
on $\Omega$ is isomorphic to the action of $G_q$ on $G_q/X$. In particular,
   $Y = \cap \, gXg^{-1} \subseteq X$.\qed

\section{Normal subgroups of $G_q$}
\subsection{The main results}Let $X$ be a subgroup of finite index of $G_q$.
The main purpose of this subsection is to show that
 $X$ is a normal subgroup of $G_q$ if and only if
  $|G_X| = [G_q : X]$.

 \smallskip
 \noindent {\bf Lemma 6.1.}  {\em
  Let $X$ be a subgroup of finite index of $G_q$.
   Then $C_{S_{\Omega}}(G_X) \cong N_{G_q}(X)/X$, Aut$\, M_X
  = C_{S_{\Omega}}(G_X)$.}

 \smallskip
 \noindent {\em Proof.}
  Since the action of
    $G_X$ on $\Omega$ is isomorphic to the action of $G_q$
     on $G_q/X$ and $G_q$ acts transitively on $G_q/X$,
     $G_X$ acts transitively on $\Omega$.
      By Lemma 2.1,
 $C_{S_{\Omega}}(G_X) \cong N_{G_X}(G_{\hat d})/G_{\hat d}$
    where $G_{\hat d}$ is the one point stabiliser of $\hat d\in \Omega$.
     Apply the fact that
     the action of
    $G_X$ on $\Omega$ is isomorphic to the action of $G_q$
     on $G_q/X$ one more time, one has
      $N_{G_X}(G_{\hat d} )/G_{\hat d} \cong N_{G_q}(X)/X$.

    \smallskip  Aut$\,M_X =  C_{S_{\Omega}}(G_X) $ follows
  from the fact that the incidence relations of $(V, E,F)$ is completely determined by $G_X$
   and that $\sigma\in S_{\Omega}$ is an automorphism if and only if $\sigma$ preserves the  incidence
    relations of $(V, E,F)$.
 \qed

   \smallskip

   \noindent {\bf Proposition 6.2.} {\em
        Let $X$ be a   subgroup of $G_q$ of finite index and
         let $G_{X}= \left <r_1, r_2\right>$.
        Then $X$ is a normal subgroup of $G_q$ if and only if
          $|G_X|= [G_q : X]$.
      In the case $X$ is normal, $G_X \cong C_{S_{\Omega}} (G_X)\cong Aut\, M_X$.
         }

  \smallskip
  \noindent {\em Proof.}  
    Suppose that  $ |G_X|= [G_q : X]$.
     By (5.2),  one has
   $$|G_q/Y |= |G_q/ker\,f|= |G_X| = [G_q : X].\eqno(6.1)$$
   By Lemma 5.4, $Y \subseteq X$. Hence  $X= Y$ is normal.
  Conversely, suppose that $X$ is a normal subgroup.
  By (5.2) and Lemma 5.4, one has $G_X \cong G_q/Y =
  G_q/X$. Further, one has $G_X \cong C_{S_{\Omega}} (G_X)$ (see Corollary 2.2).
   This completes the proof of the proposition.\qed


\subsection{ The geometric invariance of $X$
 and $M_X$} The invariance of $X$ and $M_X$
 (see (xiii) of subsection 3.1)
 can be described by the group $G_X =\left < r_0, r_1, r_2\right >=\left <r_1, r_2\right >
 $ as follows.
 \begin{enumerate}
 \item[(i)]
  $|E|=$ the number of disjoint cycles (counting 1-cycles) in $r_1$,
   $[G_q : X] = |\Omega|=$ the number of entries (counting 1-cycles) in $r_1$ and
    $\tau_2=$ the number of 1-cycles of $r_1$.
   \item[(ii)]
     $|F|= $ the number of disjoint cycles (counting 1-cycles) in $r_2$.
     For each divisor $r$ of $q$ ($1 <r \le q$), $v_r= $ the number of $q/r$
      cycles in $r_2$.
     \item[(iii)] $|V|=$ the number of disjoint cycles   (counting 1-cycles)  in $r_0
      = r_2^{-1}r_1^{-1}$.
    Each cycle in $r_0$ represents a vertex, the degree of the vertex
     is the length of the cycle.
    \item[(iv)]  $N_{G_q} (X)/X \cong  C_{S_{\Omega}}(G_X) = Aut\,M_X$.
     $X$ is normal  iff
   $ C_{S_{\Omega}}(G_X)\cong G_X\cong G_q /X$.
    \end{enumerate}


\allowdisplaybreaks

\begin{center}
\begin{figure}
\beginpicture

\setcoordinatesystem units <6pt,6pt>

\setplotarea x from  -18.5 to 63, y from -1 to 22

\circulararc -180 degrees from 0 0 center at 3.75 0
\circulararc -180 degrees from 7.5 0 center at 11.25 0

\circulararc -180 degrees from 15 0 center at 18.75 0
\circulararc -180 degrees from 22.5 0 center at 26.25 0


\setlinear \plot 0 0 30 0 /

\setlinear \plot  0   0    0 18 /
\setlinear \plot 15  0    15 18 /
\setlinear \plot 30  0    30 18 /

\put {$1$} at   2 13
\put {$4$} at   17 13
\put {$3$} at  32  13
\put {$3$} at 18.75 1.75
\put {$2$} at  26.25 1.75
\put {$1$} at 11.25 1.75
\put {$2$} at  3.75 1.75


\put {$  \vee$} at   0 8
\put {$\vee  $} at   15 8
\put {$\vee   $} at   30 8


\put {$ > $} at 18.75 3.75
\put {$< $} at  26.25 3.75
\put {$ < $} at 11.25 3.75
\put {$> $} at  3.75 3.75

\put {$\scriptsize{0}$} at 0 -3
\put {$\scriptsize{1/\lambda}$} at 7.5 -3

\put {$\scriptsize{\lambda }$} at 15 -3
\put {$\scriptsize{2\lambda}$} at 30 -3
\put {$\scriptsize{3/\lambda}$} at 22.5 -3


\put {\bf Figure 1} at 16 -7
\endpicture
\end{figure}
\end{center}

\section{ First  Application : Identifying Normal Subgroups}
\subsection{Identifying normal subgroups} Let $X$ be a subgroup of $G_q$
 of finite index and let $G_X$
  be given as in (5.1).
 The construction
 of $r_0$ is  difficult as it is not an easy matter to list the
directed   edges of $c_{v}$ according to the orientation of $M_X$.
{\em However,
it is an  easy matter to write down the permutation representations of
 $r_1$ and $r_2$} (see Examples 7.1-7.5).
  As a consequence, the group $G_X$ can be determined by $M_X$ without any membership test.
  The
  order of  $G_X =\left <  r_1, r_2 \right >    $
 can be determined  by {\bf GAP}.
Hence whether  $X$ is normal can be determined  by
  Proposition 6.2 and  {\bf GAP} [G].


\subsection{Normalisers and automorphisms}
   Let $N(X)$ be the normaliser of $X$ in $PSL(2, \Bbb R)$.
  In the case $q\ne 3,4,6$, the Margulis characterisation of arithmeticity in terms
      of the commensurator implies that
 $N_{G_q}(X)  =N(X)$. Hence $N(X)/X$  can be determined by Lemma 6.1.
The automorphism group of the map $M_X$ can be determined by Lemma 6.1 also.

\subsection{ Discussion} (i)
 A map is called {\em quasi-regular} if every vertex has the same degree, every face
     has the same number of edges, and  either it has no free edges or
 all its edges are free.
  It is clear that $X$ is normal only if $M_X$ is quasi-regular. The converse is not
      true (see Example 7.1).
(ii) A map $M$ is {\em regular} if Aut$\, M$ acts transitively on the set of darts.
A  quasi-regular map on a
genus zero surface is regular (see Corollary 6.4 of [JS]).
 (iii) By our result in subsection 4.3, $|\Omega| = [G_q: X]$.
  As a consequence, $M_X$ is regular if and only if $X$ is normal (see Lemma 6.1).


 \smallskip

 \noindent {\bf Example 7.1.} Let $\lambda = \sqrt2$ and let
$$M_X = \{  {-\infty }_{_{_{\smile}} }\  \hspace{-.37cm}_{ _{_{_{_{_{_{_{1}}}}}}} }
 {0}
_{_{_{\smile}} }\  \hspace{-.37cm}_{ _{_{_{_{_{_{_{2}}}}}}}}
1/\lambda
 _{_{_{\smile}}}\  \hspace{-.36cm} _{_{_{_{_{_{_{_{1}}}}}}}}
 \lambda
 _{_{_{\smile}}}\  \hspace{-.36cm} _{_{_{_{_{_{_{_{3}}}}}}}}
  3/\lambda \,
    _{_{_{\smile}}}\  \hspace{-.36cm} _{_{_{_{_{_{_{_{2}}}}}}}}
 {2\lambda}\,
  _{_{_{\smile}}} \  \hspace{-.36cm} _{_{_{_{_{_{_{_{3}}}}}}}}
 {\infty}
  \}
  \eqno(7.1)$$
be the Hecke-Farey symbols of $X\subseteq G_4$  of index 8 (see Figure 1).
Applying (i)-(iii) of subsection 3.1,  a set of independent generators
 is given by
$$  g_1=\left (
\begin{array}{rr}
\lambda  & 1\\
1 & \lambda \\
\end{array}
\right ) ,\,\,
 g_2=\left (
\begin{array}{rr}
7 & -2\lambda\\
2\lambda & -1 \\
\end{array}
\right ) ,\,\,
g_3=\left (
\begin{array}{rr}
3\lambda & -7\\
1 & -\lambda \\
\end{array}
\right ) ,\eqno(7.2)$$
where the side pairing $g_r$      pairs the even lines with label $r$.
One sees easily from Figure 1 that $$r_2 =
 (1,2, \bar 1, \bar 4)(  4,3, \bar 2, \bar3),\,
r_1 =
(1,\bar 1) (2,\bar 2)
(3,\bar 3) (4,\bar 4),\eqno(7.3)$$
where $x$ and $\bar x$ are the  directed edges associated with the same edge.
 It follows that
  $r_0 = r_2^{-1} r_1^{-1}= (2,3, \bar
   2,  1)(\bar 4, \bar 3, 4, \bar 1)$.
   By {\bf GAP}, one has  $|G_X| = 16 >  8  = [G_4 : X]$.
   By Proposition 6.2,  $X$ is not normal.
Note that $M_X$ is quasi-regular but not regular.  Further,  $Aut\, M_X =C_{S_{8}} (G_{X})\cong
 Z_2 \times Z_2$.
The following  map  is also  quasi-regular.
$M_Y = \{  {-\infty }_{_{_{\smile}} }\  \hspace{-.37cm}_{ _{_{_{_{_{_{_{1}}}}}}} }
 {0}
_{_{_{\smile}} }\  \hspace{-.37cm}_{ _{_{_{_{_{_{_{3}}}}}}}}
1/2
 _{_{_{\smile}}}\  \hspace{-.36cm} _{_{_{_{_{_{_{_{1}}}}}}}}
 1
 _{_{_{\smile}}}\  \hspace{-.36cm} _{_{_{_{_{_{_{_{2}}}}}}}}
  2
    _{_{_{\smile}}}\  \hspace{-.36cm} _{_{_{_{_{_{_{_{3}}}}}}}}
 {3} _{_{_{\smile}}} \  \hspace{-.36cm} _{_{_{_{_{_{_{_{2}}}}}}}}
 {\infty}
  \}.$
   Note that
   $[PSL(2, \Bbb Z): Y] = 12$ and that $Y$ is not normal.

\smallskip

 \noindent {\bf Example 7.2.} Let
$$M_X = \{  {-\infty }_{_{_{\smile}} }\  \hspace{-.37cm}_{ _{_{_{_{_{_{_{1}}}}}}} }  {0}
_{_{_{\smile}} }\  \hspace{-.37cm}_{ _{_{_{_{_{_{_{2}}}}}}}}
1
 _{_{_{\smile}}}\  \hspace{-.36cm} _{_{_{_{_{_{_{_{2}}}}}}}}
 3/2
 _{_{_{\smile}}}\  \hspace{-.36cm} _{_{_{_{_{_{_{_{3}}}}}}}}
  2
    _{_{_{\smile}}}\  \hspace{-.36cm} _{_{_{_{_{_{_{_{3}}}}}}}}
 {3} _{_{_{\smile}}} \  \hspace{-.36cm} _{_{_{_{_{_{_{_{1}}}}}}}}
 {\infty}
  \}
   \eqno(7.4)$$
be the Hecke-Farey symbols of $X \subseteq \Gamma = PSL(2, \Bbb Z)$ of index 12
(see Figure 2).
Applying (i)-(iii) of subsection 3.1,  a set of independent generators
 is given by
{ $$
  g_1=\left (
\begin{array}{rr}
1 & 3\\
0 & 1 \\
\end{array}
\right ) ,\,\,
 g_2=\left (
\begin{array}{rr}
2 & -3\\
3 & -4 \\
\end{array}
\right ) ,\,\,
g_3=\left (
\begin{array}{rr}
5 & -12\\
3 & -7 \\
\end{array}
\right ) .\eqno(7.5)$$}
One sees easily from Figure 2 that  $$r_2 =
 (1,2, \bar 4)(  4,6, \bar 5)( 5, \bar 3,\bar 1 )(\bar 6,\bar 2, 3),\,
 r_1 = (1,\bar 1) (2,\bar 2)
(3,\bar 3) (4,\bar 4)
(5,\bar 5) (6,\bar 6),\eqno(7.6)$$
where $x$ and $\bar x$ are the two directed edges associated with the same edge.
By {\bf GAP},  one has
   $  |G_X| = [\Gamma : X] = 12$.
      By Proposition 6.2, $X$ is a normal subgroup of $\Gamma$.
 Note that $X= \Gamma (3)$.

\allowdisplaybreaks

\begin{center}
\begin{figure}
\beginpicture

\setcoordinatesystem units <5pt,5pt>

\setplotarea x from  -17 to 63, y from -1 to 28

\circulararc -180 degrees from 0 0 center at 7.5 0
\circulararc -180 degrees from 15 0 center at 22.5 0
\circulararc -180 degrees from 15 0 center at 18.75 0
\circulararc -180 degrees from 22.5 0 center at 26.25 0
\circulararc -180 degrees from 30 0 center at 37.5 0

\setlinear \plot 0 0 45 0 /

\setlinear \plot  0   0 0 25 /
\setlinear \plot 15  0 15 25 /
\setlinear \plot 30  0 30 25 /
\setlinear \plot 45  0 45 25 /

\put {$1$} at   2 15
\put {$1$} at  43  15
\put {$4$} at   17 15
\put {$5$} at  32  15

\put {$2$} at  7.3 5.5
\put {$2$} at 18.75 1.75
\put {$3$} at  26.25 1.75
\put {$3$} at  37.5  5.5
\put {$6$} at  22.5  5.5



\put {$   \vee$} at   0 12
\put {$   \vee$} at   45 12
\put {$\vee   $} at   15 12
\put {$\vee   $} at   30 12

\put {${ > } $ } at  7.5 7.5
\put {$ >$} at  22.5 7.5
\put {$ <$} at  37.5 7.5
\put {$ < $} at 18.75 3.75
\put {$> $} at  26.25 3.75

\put {$\scriptsize{0}$} at 0 -3
\put {$\scriptsize{1}$} at 15 -3
\put {$\scriptsize{2}$} at 30 -3
\put {$\scriptsize{3/2}$} at 22.5 -3
\put {$\scriptsize{3}$} at 45 -3

\put {\bf Figure 2} at 22.5 -7
\endpicture
\end{figure}
\end{center}

\vspace{-.4cm}

\noindent {\bf Example 7.3.}
Let $\lambda = 2cos\pi/5 = (1+\sqrt 5)/2$ and let
 $$M_X = \{
 -\infty\,
 _{_{_{_\smile}} }\  \hspace{-.37cm}_{ _{_{_{_{_{_{_{\bullet}}}}}}} }
{0/1 \,}
 _{_{_{\smile}} }\  \hspace{-.37cm}_{ _{_{_{_{_{_{_{\,\bullet}}}}}}} }
{1/\lambda\, }
 _{_{_{\smile}} }\  \hspace{-.37cm}_{ _{_{_{_{_{_{_{\,\bullet}}}}}}} }
{\lambda/\lambda\, }
 _{_{_{\smile}} }\  \hspace{-.37cm}_{ _{_{_{_{_{_{_{\,\bullet}}}}}}} }
{\lambda/1\, }
 _{_{_{\smile}} }\  \hspace{-.37cm}_{ _{_{_{_{_{_{_{\,\bullet}}}}}}} }
{\infty}
  \}
  \eqno(7.7)
  $$
  be the Hecke-Farey symbols of $X \subseteq G_5$ of index 10.
   Similar to Examples 7.1 and 7.2 (with appropriate labeling of the even lines), we have
   $
   r_1 = (1, \bar 1) (2, \bar 2)(3, \bar{3})
   (4, \bar{4}) (5, \bar{5})$
    and
   $
   r_2= (1,2,3,4,5)(\bar1 )(\bar2 )(\bar3 )(\bar4)(\bar 5)$.
    In particular, $M_X$ consists of one 5-gon and five 1-gon.
    By subsection 7.3, $X$ is not normal.
   By
    {\bf GAP},
   $C_{S_{10}}(\left < r_1, r_2\right >) =
   \left < (1,2,3,4,5)(\bar1,\bar2,\bar3 ,\bar4,\bar 5)\right >$.
Hence
$Aut\,M_X =   C_{S_{10}}(\left < r_1, r_2\right >) \cong N(X)/X  \cong Z_5$
 (see subsection 7.2 and Lemma 6.1).

\smallskip
\noindent {\bf Example 7.4.} Throughout the example,
  $\lambda = \sqrt 3$. (i)  Let $X$ be given as in subsection 3.3.
 $X$ has index 3 in $G_6$.
  The edges of $M_X$ are free edges (see Figure 3).
   The  set of darts is $\{1,2,3\}$.
  It follows easily that
  $r_1= (1)(2)(3)$ and that $r_2 = (1,2, 3)$. As a consequence, $G_X =
   C_{S_3}(G_{X}) = \left < (1,2,3)\right >$.
   By Proposition 6.2,
     $X$ is a normal subgroup of $G_6$.
 (ii)
Let
$ M_Y= \{
   {-\infty} \,\,_{_{_{_{\smile}} }}\  \hspace{-.37cm}_{ _{_{_{_{_{_{_{1}}}}}}} }
   {0} \,\,_{_{_{_{\smile}} }}\  \hspace{-.37cm}_{ _{_{_{_{_{_{_{\circ}}}}}}} }
 {1/\lambda} \,\,_{_{_{_{\smile}} }}\  \hspace{-.37cm}_{ _{_{_{_{_{_{_{1}}}}}}} }
     {\lambda/2 } \,\,_{_{_{_{\smile}} }}\  \hspace{-.37cm}_{ _{_{_{_{_{_{_{e_2}}}}}}} }
  \infty   %
\}$ (see Figure 4).
 $Y$ is a subgroup of index 3 in $G_6$.
 An easy calculation
      of the side pairings shows that $Y = G_0(2)=\{
       (x_{ij} )\in G_6 : x_{12} \equiv 0\,(mod\,\,2)\}$.
        By Proposition 6.2, $Y$ is not normal. Note that
         both $M_X$ and $M_Y$ are 3-gons (with different
          side pairings).

\allowdisplaybreaks

\begin{center}
\begin{figure}
\beginpicture

\setcoordinatesystem units <6pt,6pt>
\setplotarea x from  -11 to 63, y from -1 to 22

\circulararc -83 degrees from 0 0 center at 15 0
\circulararc -85 degrees from 0 0 center at 3.75 0

\circulararc -78 degrees from 4.00 3.75  center at 3.75 0

\circulararc -41 degrees from 7.5 0 center at 30  0

\circulararc -85 degrees from 7.5 0 center at 10.35 0
\circulararc -80 degrees from 10.65 2.85  center at 10.35 0


\setlinear \plot  0   15.4    0 20 /
\setlinear \plot  0   0       0 14.7 /

\setlinear \plot 13.2  0  13.2 20 /

\put {$0$} at   0 -3
\put {$1/\lambda$} at   7.5 -3
\put {$\lambda/2$} at   13.2 -3

\put {$ e_2$} at   15 15

\put {$\circ$} at   0 15
\put {$\circ$} at   3.75 3.75
\put {$\circ$} at   10.35 2.85

\put {$1$} at   1.5 15
\put {$2$} at   3.75 5.5
\put {$3$} at   10.35 4.5


\circulararc -83 degrees from 30 0 center at 45 0
\circulararc -85 degrees from 30 0 center at 33.75 0

\circulararc -80 degrees from 34.00 3.75  center at 33.75 0

\circulararc -41 degrees from 37.5 0 center at 60  0

\circulararc -180 degrees from 37.5 0 center at 40.35 0

\setlinear \plot  30   0    30 20 /

\setlinear \plot 43.2  0  43.2 20 /

\put {$0$} at   30 -3
\put {$1/\lambda$} at   37.5 -3
\put {$\lambda/2$} at   43.2 -3

\put {$ e_2$} at   45 15

\put {$\vee$} at   30 15
\put {$\circ$} at   33.75 3.75
\put {$<$} at   40.35 2.85

\put {$1$} at   31.5 15
\put {$2$} at   33.75 5.5
\put {$1$} at   40.35 4.5











\put {\bf Figure 3} at 7.5 -7
\put {\bf Figure 4} at 37.5 -7
\endpicture
\end{figure}
\end{center}

\vspace{-.4cm}



\noindent {\bf Example 7.5.  Normal non-congruence subgroups of index 42 of $PSL(2, \Bbb Z)$.}
 Let $X$  be the group associated with   $ M_{1,6}$
  and let $Y$ be the normal map group of $X$ (see (3.10) and Definition 5.1).
 By {\bf GAP}, the group
 $G_X$  has order 42.
 Hence $Y$ is a normal subgroup of index 42.
  Since $X$ is non-congruence,
  $Y$ is also non-congruence.
   Note that 42  is  the smallest possible index of a normal
    non-congruence subgroup of $PSL(2, \Bbb Z)$ (see pp. 276 of [N]).


\section{Second Application :  Congruence Subgroup Problem for $PSL(2, \Bbb Z)$}
 The {\em principal congruence subgroup} of $G_q$ associated with  $(r)\subseteq \Bbb Z[\lambda_q]$ is the subgroup
   $G(r) = \{ x\in G_q \,:\, x \equiv \pm I\,\,(mod\,\,  r)\}$.
 $X\subseteq G_q$ is a {\em congruence subgroup} if $G(r) \subseteq X$ for some $r$.
Let $X$ be a subgroup of finite index of $G_q$ and let $n$ be the least common multiple
 of the degrees of the vertices. The number $n$ is called the {\em level} of $X$.
    It is known that if $X\subseteq G_q$ has level $n$, then
    \begin{enumerate}
    \item[(i)] if $q = 3,4,6$, then
     $X$ is congruence if and only if $ G(n\lambda_q) \subseteq X$
      (see [Wo], [P]),
      \item[(ii)] if $q=5$, then
    $X$ is congruence if and only if $ G(2n) \subseteq X$
     (see [LL]).
     \end{enumerate}
      As a consequence, the algorithm  in [LLT3] that
      determines whether a group is congruence extends essentially
      verbatim to the Hecke groups $G_q$ when $q \le 6$
       (see pp. 46 of [HR]). The weakness  of their
       algorithm  is that the calculation becomes very lengthy when the level $n$ is
        large. In the case $q = 3$,  to the best of our knowledge, the congruence test
         for $G_3 = \Gamma = PSL(2, \Bbb Z)$ developed by
         Tim Hsu [H] is the most effective one in the literature. His test
          can be implemented easily as long as  a $T U$-representation
           of $X$ is determined,
           where
           $$T =\left (
\begin{array}{rr}
1 & 1\\
0 &1 \\
\end{array}
\right ) ,\,
 U =\left (
\begin{array}{rr}
1 & 0\\
1 &1\\
\end{array}
\right ).\eqno(8.1)$$
Recall that a $TU$-representation is a set of two permutations \{$f(T), f(U)\}$,
 where $f(x)$ denotes the permutation representation of $x\in\{ T,U\}$ on $PSL(2, \Bbb Z)/X$.
The determination of $f(T)$ and $f(U)$  involves matrix
  multiplication and membership test.
   As this procedure may be tedious when the index
      $[PSL(2,\Bbb Z) :X]$
      is large (see Section 5.4 of  [KL] for example), we suggest the following,
      which makes the determination of whether $X$ is congruence very easy
       as long as $X$ is given as in subsection 3.1.

 \smallskip
\noindent {\bf Proposition 8.1.} {\em Let  $G_X = \left <r_0, r_1, r_2\right >$
   be the monodromy group  of $X\subseteq PSL(2, \Bbb Z)$. Then
      $$f(T)  = r_0 = r_2^{-1} r_1^{-1} ,\,\, f(U) = r_1 r_0^{-1} r_1^{-1}.
      \eqno(8.2)$$
       $X$ is a congruence subgroup  if and only if
        of $f(T)$ and $f(U)$ satisfy the conditions
        given in Section $3$ of $[H]$.}

\smallskip
\noindent {\em Proof.} By (5.2) and Lemma 5.4, $r_1$ gives the action of
{\tiny $\left (
\begin{array}{rr}
0 & 1\\
-1 & 0\\
\end{array}
\right )$}  and $r_2$ gives the action of {\tiny $\left (
\begin{array}{rr}
0 & 1\\
-1 & 1 \\
\end{array}
\right )$} on $PSL(2, \Bbb Z)/X$. It follows by direct calculation that
 $f(T) = r_0$ and that $f(U) = r_1 r_0^{-1} r_1^{-1}$. The proposition now follows
  by applying Theorem 3.1 of [H].\qed

 \smallskip
To ensure the reader that our procedure described as above is practical
 and makes Hsu's algorithm a perfect choice when $X$ is given as in subsection 3.1, we
  provide
 the following.

\smallskip
\noindent {\bf Example 8.2.}
 Let
$$ M_X = \{  {-\infty }_{_{_{\smile}} }\  \hspace{-.37cm}_{ _{_{_{_{_{_{_{\bullet}}}}}}} }
 0
_{_{_{\smile}} }\  \hspace{-.37cm}_{ _{_{_{_{_{_{_{\circ}}}}}}}}
1
 _{_{_{\smile}}}\  \hspace{-.36cm} _{_{_{_{_{_{_{_{\circ}}}}}}}}
 2
_{_{_{\smile}} }\  \hspace{-.37cm}_{ _{_{_{_{_{_{_{\circ}}}}}}}}
3
 _{_{_{\smile}}}\  \hspace{-.36cm} _{_{_{_{_{_{_{_{\bullet}}}}}}}}
 \infty
  \}\eqno(8.3)$$
be the Hecke-Farey symbols of a subgroup $X$ of index 11 of $PSL(2, \Bbb Z)$
 (see Figure 5).
 Then
 $$ r_1= (1, \bar 1)(2,\bar 2)(3, \bar 3)(4, \bar 4),\,\,
  r_2 = (\bar 1) (1,5,\bar 2)(2,6,\bar 3)(3,7,\bar 4)( 4).\eqno(8.4)$$
  Apply Proposition A.1, the $TU$-representation of $X$ is given as follows.
   $$ f( T)= (\bar 1, \bar 2, \bar 3, \bar 4, 4,7,3,6,2,5,1),\,\,
  f(U) = (\bar 1,5,\bar 2, 6, \bar 3, 7, \bar 4, 4,3,2,1)
  .\eqno(8.5)$$
 By the algorithm given in Section 3 of [H], $X$ is non-congruence.
 By {\bf GAP},  $\left
  < f(T),f( U)\right > \cong A_{11}$ the alternating group on 11 letters.
   The  group $\cap \,gXg^{-1} $  was first studied by Magnus [M]
    as part of his study of non-congruence subgroups of $PSL(2, \Bbb Z)$.

\allowdisplaybreaks

\begin{center}
\begin{figure}
\beginpicture

\setcoordinatesystem units <4pt,4pt>


\setplotarea x from   -26.5 to 83, y from -1 to 28

\circulararc -84 degrees from 0 0 center at 7.5 0

\circulararc  -90 degrees from 7.5 7.5 center at 7.5 0

\circulararc -87 degrees from 15 0 center at 22.5 0
\circulararc  -88 degrees from 22.9 7.5 center at 22.5 0


\circulararc -87 degrees from 30 0 center at 37.5 0
\circulararc -87 degrees from 38 7.5 center at 37.5 0


\circulararc -60 degrees from 45 0 center at  60 0
\circulararc  60 degrees from 0 0 center at  -15 0

\setlinear \plot 0 0 45 0 /

\setlinear \plot  0   0 0 25 /
\setlinear \plot 15  0 15 25 /
\setlinear \plot 30  0 30 25 /
\setlinear \plot 45  0 45 25 /


\setlinear \plot -7.5  13 -7.5 25 /
\setlinear \plot 52.5  13  52.5 25 /

\put {$\bullet$} at   -7.5 13
\put {$\bullet$} at  52.5  13

\put {$\vee$} at   0 13
\put {$\vee$} at  45  13

\put {$2$} at   17 13
\put {$3$} at  32  13

\put {$1$} at   2  13
\put {$4$} at  43  13

\put {$5$} at  7.3 5.5
\put {$7$} at  37.5  5.5
\put {$6$} at  22.5  5.5



\put {$\vee   $} at   15 13
\put {$\vee   $} at   30 13

\put {${ \circ } $ } at  7.5 7.5
\put {$ \circ$} at  22.5 7.5
\put {$ \circ$} at  37.5 7.5

\put {$\scriptsize{0}$} at 0 -3
\put {$\scriptsize{1}$} at 15 -3
\put {$\scriptsize{2}$} at 30 -3
\put {$\scriptsize{3}$} at 45 -3

\put {\bf Figure 5} at 22.5 -7
\endpicture
\end{figure}
\end{center}

\vspace{-.5cm}

\section*{Appendix A}

The following lemma is   an application of our
 study of the Hecke-Farey symbols of the Hecke groups.
 To the best of our knowledge, this lemma  was first proved in [CS]
 with the help of the celebrated Selberg's Theorem. Our proof is more
  geometrical and does not make use of Selberg's Theorem.

\smallskip
\noindent {\bf Lemma A1.} {\em Let $q$ be composite. Then
 $G_q$ has infinitely many normal subgroups of finite indices with
  torsion.}

\smallskip
\noindent {\em Proof.}
Recall that $G_q$ is a free product of $Z_2 =\left <S\right >$ and  $Z_q =\left < R_q\right >$, where $R _q = ST^{-1}$. Suppose that $q$ is even. Let $D_{2m}= \left
 < a, b\,:\, a^2 = b^2 = (ab)^m =1\right>$.
  Define $f_m \, :\, G_q \to D_{2m}$ by $f_m(S) = a$, $f_m(R_q)\to b$.
  It follows easily that the kernel of $f_m$ is a normal subgroup of index $2m$
   for every  $m$
   with torsion.
    In the case $q$ is odd, let $r$ be a prime divisor of $q$ and
   let $X_{2m}$
   be a subgroup of $G_r$ of index $2rm$ whose
  additional structure of
     the Hecke-Farey symbols
     takes the form
${x_i} _{_{_{\smile}} }\  \hspace{-.37cm} _{_{_{_{_{_{_{a}}}}}}}  x_{i+1}
 $ only (a special polygon of $X_{2m}$ is a union of $2m$ $r$-gons with $2m(r-2)+2$ even lines).
  It follows that $X_{2m}$ is torsion free
   (see  subsection 3.1). As a consequence, $ K_{2m} =\cap _{g\in G_r} g
   X_{2m}g^{-1}$ is  normal torsion free of finite index of $G_r$
   and $G_r/K_{2m}$ is a finite  group generated by $S$ and $R_r$ of orders 2 and $r$
    respectively
    modulo $K_{2m}$.
    Define $f\,:\, G_q\to G_r/K_{2m}$ by $f(S) = S, f(R_q) = R_r.$
  It is clear that $f$ is a surjective homomorphism and that
   $G_q$ possesses a normal subgroup $V_{2m}$  such that
    $G_q/V_{2m} \cong G_r/K_{2m}$.
    Note that $(R_q)^r \in V_{2m}$  for all $m$.
     In particular, $V_{2m}$ is torsion.
   Consequently, $G_q$ has infinitely many normal subgroups of finite indices with torsion.
       \qed

 \smallskip
   \noindent {\bf Lemma A2.}
    {\em Let $q\ge 3$ be a prime and $N$  a normal subgroup of finite index of $G_q$. Then
    $N$ is either
        $G_q$,  $G_q^2$ or  $G_q^q$.}

\smallskip
\noindent {\em Proof.}
Suppose that $N\triangleleft  \,G_q$ has
 an elliptic element $g$. It is well known  that $\left <g\right >
 $ is conjugate to either $\left <S\right >
 $  or $\left <ST^{-1}\right >
 $.
   Since $N$ is normal,
   $N$ contains all the conjugates of $\left <g\right >
 $. Hence
   $  G_q^q = \left < xSx^{-1} \,:\, x \in G_q\right >  \subseteq N$ or
  $ G_q^2 =\left < x(ST^{-1})x^{-1} \,:\, x \in G_q\right >  \subseteq N$.
   It is easy to see that $G_q^q$ is a  normal subgroup of index $q$ whose special polygon is the $q$-gon
 $$ \{  {-\infty }_{_{_{\smile}} }\  \hspace{-.37cm}_{ _{_{_{_{_{_{_{\circ}}}}}}} }
 {0}
_{_{_{\smile}} }\  \hspace{-.37cm}_{ _{_{_{_{_{_{_{\circ}}}}}}}}
{1/\lambda_q}
 _{_{_{\smile}}}\  \hspace{-.36cm} _{_{_{_{_{_{_{_{\circ}}}}}}}}  \cdots \cdots \,
   _{_{_{\smile}}}\        \hspace{-.37cm} _{_{_{_{_{_{_{_{\circ}}}}}}}}
 {\lambda_q/1}
_{_{_{\smile}} }\  \hspace{-.37cm}_{ _{_{_{_{_{_{_{\circ}}}}}}}}
 {\infty}
 \} \eqno(A1)$$
       and  that $G_q^2$
     is a normal subgroup of index 2 whose special polygon is the union of two
      special triangles
      $ \{  {-\infty }_{_{_{\smile}} }\  \hspace{-.37cm}_{ _{_{_{_{_{_{_{\bullet}}}}}}} }
 {0}
_{_{_{\smile}} }\  \hspace{-.37cm}_{ _{_{_{_{_{_{_{\bullet}}}}}}}}
\infty
\}.$
It follows that   the only torsion normal subgroups of $G_q$ are $G_q$, $G_q^2$ and
 $G_q^q$.\qed

\bigskip{\small

\noindent Cheng Lien Lang\\
\noindent Department of Mathematics, I-Shou  University, Kaohsiung, Taiwan.

\noindent   \texttt{cllang@isu.edu.tw}

\smallskip
\noindent Mong Lung Lang \\
\noindent Singapore 669608,
Singapore.

\noindent \texttt{lang2to46@gmail.com}}

\medskip


\end{document}